\normalbaselineskip=1.6\normalbaselineskip\normalbaselines
\magnification=1200

\def\pmb#1{\setbox0=\hbox{#1}%
 \kern-.025em\copy0\kern-\wd0
 \kern.05em\copy0\kern-\wd0
 \kern-.025em\raise.0433em\box0 }

\def\AA{{\bf A}}
\def\Gm{{\bf G}_m}
\def\Ga{{\bf G}_a}

\def\N{{\bf N}}
\def\Z{{\bf Z}}
 
\def\Q{{\bf Q}}

\def\F{{\bf F}}

\def\isom{\cong}
\def\Aut{\mathop{\hbox{Aut}}}

\def\mod{\mathop{\rm mod}\nolimits}

\def\Spec{\mathop{\rm Spec}\nolimits}
\def \bs{\bigskip}
\def\Gal{\mathop{\rm Gal}\nolimits}
\def\End{\mathop{\rm End}\nolimits}
\def\Fix{\mathop{\rm Fix}\nolimits}
\def\Aut{\mathop{\rm Aut}\nolimits}

\def \bs{\bigskip}

\def \s{\sigma}
\def \a{\alpha}

\def \G{\Gamma}

\def \con{\equiv }

\def \bF{{\bar{\bf F}}_p}

\centerline{\bf Difference subgroups of commutative
algebraic groups over finite fields }
\medskip
\centerline{\bf Thomas Scanlon and Jos\'e Felipe Voloch}
\bs

{\bf Introduction}

\bs 

The work of Chatzidakis and Hrushovski [CH] on the model theory of 
difference fields in characteristic zero showed that 
groups defined by difference equations
have a very restricted structure. For instance, if $G$ is a semi-abelian 
variety over a difference field of characteristic zero 
and $\G \subset G$ is a subgroup of ``modular type'', then for
any subvariety $X \subset G$, $X \cap \G$ is a finite union of cosets of
subgroups of $\G $ (see [C]).
Using such facts once can resolve diophantine questions
about special subgroups of $G$ (for instance, the torsion
subgroup [Hr]).
Recent work of Chatzidakis, Hrushovski and Peterzil [CHP] extends the class of
difference fields for which this sort of result is known to positive 
characteristic. In this note, we analyze the subgroups of
the torsion points of simple commutative algebraic groups 
over finite fields that can be
constructed by such difference equations. Our results are reasonably 
complete modulo some well-known conjectures in Number Theory. 
In one case, we need 
the $p$-adic version of the four exponentials conjecture and in another
we need a version of Artin's conjecture on primitive roots.

We recover part of a theorem of Boxall [Bo]
on the intersection of varieties
with the group of $m$-power torsion points, but in general this theorem
does {\it not} follow from the model-theoretic analysis, because there 
may be no field automorphism $\s$ so that the $m$-power torsion group
is contained in a modular group definable with $\s$.
On the other hand, some of the groups defined by modular difference 
equations are much larger than the group of $m$-power torsion points, so
our results are stronger in another direction. In some ways, the model 
theoretic approach extends the approach of Bogomolov 
and the original one of Lang.

\bs

{\bf 1. General set up}

\bs

A difference ring $(R,\s)$ is a commutative ring $R$ with unity,
together with an endomorphism $\s$ of $R$. We define 
$\Fix(\s)=\{ x \in R \mid \s(x)=x \}$, the subring of $R$ fixed by $\s$.
If $X$ is a scheme over $R$ we define $X^{\s}$ to be $X \times_{\Spec R}
\Spec R$, where we take $\s^* : \Spec R \to \Spec R$ to form the fiber 
product. Then $\s$ induces a map on points $X(R) \to X^{\s}(R)$ again
denoted by $\s$.  When $X \subseteq \AA^N$ is given as a subscheme of 
affine space, then the map $X(R) \to X^\s (R)$ is induced by the 
map $(x_1, \ldots, x_N) \mapsto (\sigma(x_1), \ldots, \sigma(x_N))$.
If $X$ is a group scheme then this map is a group
homomorphism. If $X$ is defined over $\Fix(\s)$ then we can identify
$X^{\s}$ with $X$.

Fix  a difference ring $(R, \s)$.
When $G$ is a commutative
group scheme  over $\Fix (\s)$ we let $\End_{\s}(G)$ denote the 
subring of the endomorphism ring of $G(R)$ 
generated by $\s:G\to G$ and $\End_R (G)$, the ring of algebraic 
group endomorphisms of $G$ over $R$. 

There
is a homomorphism $\Z[T] \to \End_{\s}(G)$ sending $T$ to $\s$.
We denote by $P(\s)$ the image of $P(T) \in \Z[T]$ by this map.

Let $p$ be a prime number. A $p$-Weil number is an algebraic number
$\a$ all of whose archimedian absolute values are all equal to
$p^r$ for some rational
number $r$. We usually omit the $p$, when it is clear from the context.
A polynomial $P(T) \in \Z[T]$ will be called {\it modular} if none of its
roots is a Weil number. The name modular is chosen to conform with the
results of [CHP] asserting that when $P(T) \in \Z[T]$ is modular, then
if $(K, \s)$ is an existentially closed difference field and $G$ 
is a semi-abelian variety over $\Fix(\s)$, then  $\ker P(\s)$
is a group of modular type.  

The reader should consult [C] for 
a more thorough discussion of modularity.  An important feature of 
modular groups is that if $\G$ is a modular subgroup of 
$G(K)$ for some algebraic group $G$ over $K$  and $X \subseteq G$
is a subvariety, then $X \cap \G$ is a finite union of cosets of 
subgroups of $\G$.

If $P(T)$ has a Weil number as root, then 
$\ker P(\s)$ contains an infinite group commensurable with a group of
the form $G(\Fix(\s^nF^m))$, for some $n,m$ integers, where $F$ is the
Frobenius endomorphism $x \mapsto x^p$.  The results of [CH] show
that every commutative subgroup of a commutative algebraic 
group defined over a finite field definable in an existentially closed
difference field is commensurable with a group of the form 
$\ker \Lambda$ where $\Lambda \in \End_{\s} (G)$.  If $G$ is a semi-abelian
variety, then the groups $\ker \Lambda$ are contained in $\ker P(\s)$
for some non-zero $P(T) \in \Z[T]$.

Our main goal is to describe the groups of the form $\ker \Lambda 
(\bar{k}, \s)$,
where $G$ is a commutative algebraic group over a finite field $k$
of characteristic $p$, $\s$ is an element of the absolute Galois group
of $k$, $\Lambda \in \End_\s (G)$, and 
$\ker \Lambda$ is of modular type over an existentially closed difference
field. 
A slightly lower target is to characterize which ``natural'' 
subgroups of $G({\bar k})$ are contained in some modular group. Here,
``natural'' means defined by some simple group-theoretic condition. In 
particular we will deal with the $l$-primary torsion of $G({\bar k})$
for any prime $l$, or more generally, groups of the form
$G_S = \{ x \in G({\bar k}) \mid mx=0, m=\prod_{p \in S}p^{c_p}\}$,
where $S$ is some set of prime numbers.
When $S$ consists of a single prime, our results are complete.

If $G$ is an algebraic group over a finite field $k$, then $G({\bar k})$
is locally finite. So it is impossible that $G({\bar k})= \G_1\cdot\ldots\cdot
\G_n$, with each $\G_i$ modular, because $G({\bar k})$ contains 
many subvarieties
other than groups subvarieties while $\G_1\cdot\ldots\cdot
\G_n$ would be modular if each $\G_i$ were. 
This puts a restriction on $S$. In particular,
if $S \cup S'$ is the set of all primes then one of $G_S,G_{S'}$ will
not be contained in a modular group.

\bs

{\bf 2. The multiplicative group}

\bs

In analyzing the case of the multiplicative group $\Gm(\bF)$ we will see that
finding appropriate equations and automorphisms comes down to solving
exponential equations $\ell$-adically. 
These equations are easy to solve in the 
case of $\Gm$, but as the rank of the endomorphism group of the algebraic 
group grows, so does the difficulty in solving these equations.

\proclaim Theorem 1. For any prime $\ell \ne p$, there exists
$\s \in \Gal(\bF/\F_p)$
and $P(T) \in \Z[T]$ modular such that the $\ell$-primary torsion subgroup
of $\Gm$ is contained in $\ker P(\s)$.

{\it Proof:} Consider the action of $\Gal(\bF/\F_p)$ on 
$T_\ell \Gm$.
The image of $\Gal(\bF/\F_p)$ on $\Aut(T_\ell \Gm) \isom \Z_\ell^\times$
is the $\ell$-adic
closure of the group generated by $p$. If $a$ is an integer
in the $\ell$-adic closure of the group generated by $p$, is not a Weil number
(that is, is not an integral power of $p$),
we take $\s$ in the Galois group mapping to $a$ and $P(T)=T-a$
and we will be done. Thus it is enough to show that such
$a$ exists. This is immediate. Since 
the $\ell$-adic closure of the group generated by $p$ is open, we can take
$a$ sufficiently close to $1$ $\ell$-adically, not a power of $p$ and
it will serve our purpose.  $\diamondsuit$

Let us investigate now for which sets $S$ of primes, $G_S$ is contained
in $\ker \s -a$, for suitable $\s \in \Gal(\bF/\F_p)$ 
and $a \in \Z_+ \setminus p^\N$.
As in the proof of the
theorem, for $\ell \in S$ we need to know if $a$ is in
the $\ell$-adic closure of the group generated by $p$ in $\Z_\ell$.
According to Artin's conjecture on primitive roots, $p$ should be
a primitive root $\mod \ell$ for a positive proportion of primes $\ell$.
Also, one expects that $p^{\ell-1} \not\con 1 (\mod \ell^2)$ for most
primes $\ell$, therefore one expects that, for a positive proportion
of primes $\ell$
the $\ell$-adic closure of the group generated by $p$ in $\Z_\ell$ is
$\Z_\ell^\times$ and thus $S$ can be taken to be a set of primes of positive
density. Unfortunately, to prove this seems beyond the power of
current techniques.

\bs
{\bf 3. Simple Abelian Varieties}
\bs

If $A$ is a simple abelian variety over a finite field $k$ then $\End A$ 
contains $\Z [F]$ as a subring where $F: A \to A$ is the 
Frobenius relative to $k$.  
For each rational prime $\ell$, we let $T_\ell A$ denote the Tate
module of $A$ and $\rho_\ell : \Gal (\bar{k} / k) 
\to \Aut_{\Z_\ell} (T_\ell A)$ the Galois action on the Tate module.
Let $\iota_\ell : \End(A) \otimes \Z_\ell \to \End_{\Z_\ell} (T_\ell A)$
be the natural inclusion.  When $\ell \neq p$, $\iota_\ell$ is 
an injection.  If we let $f \in \Gal(\bar{k} / k)$ denote the 
Frobenius, then $\rho_\ell (f) = \iota_\ell (F \otimes 1)$.   We will
simply denote this common image by $\bar{f}$.  Let $K_\ell$ be the
normal closure over $\Q_\ell$ of the field $\Q_\ell (\bar{f})$ 
considered as a subalgebra of $\End_{\Z_\ell} (T_\ell A) \otimes_{\Z_\ell}
\Q_\ell$.    We define the multiplicative rank of $A$ at $\ell$ to be
the rank of the multiplicative subgroup of $K_\ell^\times$ generated by the 
conjugates of $\bar{f}$ over $\Q_\ell$. 

\proclaim Theorem 2. Let $A$ be a simple abelian variety defined over a
finite field $k$ of characteristic $p$. Let $\ell$ be a rational prime and
let $r$ be the multiplicative rank of $A$ at $\ell$.
If $r \leq 1$ then
the $\ell$-primary torsion of $A(\bar{k})$ is contained in 
$\ker P(\s)$, for some $\s \in \Gal(\bF/\F_p)$
and $P(T) \in \Z[T]$ modular. If $r>2$, then 
the $\ell$-primary torsion of $A(\bar{k})$
is not contained in any such group and if
$r=2$ the same happens assuming the $\ell$-adic four-exponentials conjecture.

{\it Proof:} The case of $r = 0$ is trivial.
The case $r=1$ is the same as the proof of Theorem 1, replacing
$p$ there by $\bar{f}$. Assume now that $r \geq 2$.  
If there were some $\s \in \Gal (\bar{k} / k)$ and modular
$P(T) \in \Z [T]$ with $\ker P(\s) \supseteq A (\bar{k}) [\ell^\infty]$,
then $\beta :=\rho_\ell (\s)$ would have all algebraic eigenvalues on 
$T_\ell A$ and so would be an algebraic number  when considered 
as an element of $K_\ell$.  
Write $\s = f^a$ where $a \in \hat{\Z}$.  Then $\beta = \bar{f}^{a_\ell}$
where $a_\ell$ is the image of $a$ under the natural map 
$\hat{\Z} \to \Z_\ell$.    
Replacing $f$ with $f^N$ for an appropriate integer $N$, we may assume that
$\bar{f}$ is $\ell$-adically closed to $1$.  Let $\tau_1, \ldots, 
\tau_r$ be elements of $\Gal (K_\ell / \Q_\ell)$ 
such that $\tau_1 (\bar{f}), \ldots,
\tau_r(\bar{f})$ are multiplicatively independent.  
Let $x_i := \log (\tau_i (\bar{f}))$.  Let $y_1 = 1$ and
$y_2 = a_\ell$.  If $a_\ell \notin \Q$, then $y_1$ and $y_2$ are linearly
independent over $\Q$ and by the definition of $r$, the $x_i$'s
are linearly independent over $\Q$. When $r > 2$, the 
$\ell$-adic six exponential theorem asserts that at least one of 
$\exp (x_i y_j)$ is transcendental [La].  When $r = 2$, this assertion 
is the $\ell$-adic four exponentials conjecture.
Note that $\exp (x_i y_1) = \tau_i (\bar{f})$ is algebraic and 
$\exp (x_i y_2) = \tau_i (\beta)$ is also algebraic, so if $\s$ is
to satisfy any nontrivial integral equation on the $\ell$-primary
torsion, it must be a rational power of the Frobenius, but
then its eigenvalues must be $p$-Weil numbers, so it cannot 
satisfy a modular equation.  $\diamondsuit$
 
\bs

{\bf 4. The additive group}

\bs 

While for semiabelian varieties, the torsion groups 
are all contained in nontrivial 
$m$-torsion groups for some rational integer $m$,  on the 
additive group, $\Ga (\bar{k}) [m] = 0 $ or $\Ga (\bar{k})$
depending whether $(p, m) = 1$ or not.  However, $\End_k \Ga$
is a very large ring being the ring of twisted polynomials over
$k$ in one variable: $k \{ F \} := \{ \sum_{i=0}^n a_i F : 
a_i \in k, n \in \N \}$ with the commutation rule $Fa = a^p F$.
Let $\Lambda \in \End_k \Ga$.  By the $\Lambda$-power torsion 
we mean the group $\Ga (\bar{k}) [\Lambda^\infty] := 
\{ x \in \Ga (\bar{k}) : \Lambda^n (x) = 0 $ for some $n \in 
\Z_+ \}$.  
  
\proclaim Theorem 3.  Let $k$ be a finite field and $\Lambda 
\in \End_k \Ga$ be an endomorphism of $\Ga$ over $k$ 
with nontrivial kernel over $\bar{k}$.  There is no field 
automorphism $\s \in \Gal (\bar{k}/k)$ and modular group defined 
by a difference equation involving $\s$ containing the
$\Lambda$-power torsion.

{\it Proof:}  It suffices to consider the case that $\Lambda$
is irreducible over $\bar{k}$.  The $\Lambda$-Tate
module is a finite rank free  $k[[\Lambda]]$-module and its continuous
automorphism group is an order in a possibly skew finite 
extension of $k((\Lambda))$.  Let $\bar{f}$ denote the image
of the Frobenius in $\End_{cont} T_\Lambda \Ga$.  If $\s \in 
\Gal (\bar{k} / k)$ satisfies any non-trivial equation  over $\End_k \Ga$ 
on $\Ga (\bar{k}) [\Lambda^\infty]$, then the image $\beta$ of $\s$ in 
$\End_{cont} T_\Lambda \Ga \otimes_{k [[\Lambda]]}
k((\Lambda))$ would be algebraic over $\bar{f}$ and 
hence over $k(\Lambda)$.  Write $\s = f^a$ with $a \in \hat{\Z}$.
Then,  $\beta = {\bar f}^{a_p}$.   If $a \in \Q$, then $\s$ cannot satisfy 
a modular equation as on the maximal $p$ extension of $\F_p$ 
some power of $\s$ would agree with a Frobenius and groups defined 
with a Frobenius are definable in the field language so are either
finite or all of the additive group,
but if $a \notin \Q$, then by a theorem of 
Mend\`{e}s France and van der Poorten [MFvdP], $\beta$ is transcendental.    
$\diamondsuit$

\bs 

{\bf 5. Semiabelian varieties over the ring of Witt vectors}

\bs

We wish to point out some consequences of our methods and results
to semiabelian varieties over the ring of Witt vectors. Before
doing that, we need to remain in the situation we had before and
study difference equations on the additive group $\Ga$ over a finite
field. Consider
$\G = \ker P(\s)\subset \Ga(\bF)$, for some $\s \in \Gal(\bF/\F_p)$
and $P(T) \in \F_p[T]$. We may assume that $P(0) \ne 0$ and 
that $P(T) |(T^n-1)$ for some $n$. Then, $\G \subset \Fix (\s^n)$,
and $\G$ is a vector space over $\Fix(\s)$.
We thus conclude that $\G$ is finite
if $\Fix(\s)$ is finite and is either $0$ or infinite, otherwise.

Now let $W(\bF)$ be the ring of infinite Witt vectors over $\bF$ and
$G/W(\bF)$ a semiabelian scheme. If $\s \in \Gal(\bF/\F_p)$ then
$\s$ acts on $W(\bF)$ and we consider $R=(W(\bF),\s)$ as our difference
ring. Consider $\G = \ker P(\s)\subset G (W(\bF))$, for some $P(T) \in \Z[T]$.
The exact sequence $0 \to {\hat G} (pW(\bF)) 
 \to G (W(\bF)) \to {\bar G} (\bF) \to 0$,
where ${\hat G}$ is the formal group of $G$ and ${\bar G}$ is the
special fibre of $G$, reduces the study of the properties of $\G$ to 
those of the image of $\G, {\bar \G}$ in ${\bar G}$ and of
${\hat \G} = \G \cap {\hat G}$. The group ${\bar \G}$ is a group
defined by a difference equation in a semiabelian variety over a finite
field and thus the results of the previous sections apply.
The group ${\hat \G}$ is filtered by groups defined by difference equations
on additive groups over a finite field, so the above remarks apply.
In particular we get

\proclaim Theorem 4. Notation as above. If the fixed field of $\s$
in $\bF$ is finite then ${\hat \G}$ is a finitely generated $\Z_p$ module.
If $P(T)$ reduces to a monomial modulo $p$, then ${\hat \G}=0$ and $\G$
is discrete.

This follows easily from the above discussion. We point out the
case of ${\hat \G}=0$ since the methods of [Sc] then yield that
for any subvariety $X$ of $G$, there exists $c>0$ such that
for any $g \in \G$, either $g \in X$ or the $p$-adic distance from
$g$ to $X$ is bounded below by $c$.

\bs

{\bf Acknowledgements:}
The authors would like to thank MSRI for the hospitality during
the period that this work was carried out.  The authors thank
B. Conrad for his detailed comments on an earlier version of this 
paper.
The first author is supported by an NSF Postdoctoral Fellowship.
The second author would like to thank
the University of Texas's URI for financial support.

\bigskip
\centerline{\bf References.}

\noindent
[Bo] J. Boxall, {\it Sous-vari\'{e}t\'{e}s alg\'{e}briques de 
	vari\'{e}t\'{e}s semi-ab\'{e}liennes sur un corps fini},
	Number theory (Paris, 1992-1993), pp. 69--80, LMS 
	Lecture Note Ser., 215, Cambridge Univ. Press, Cambridge,
	1995.
\medskip

\noindent
[C] Z. Chatzidakis, {\it Groups definable in ACFA}, Algebraic
	Model Theory (Toronto, 1996), pp 25 -- 52, Kluwer Acad. 
	Publ., Dodrecht, 1997.
\medskip

\noindent
[CH] Z. Chatzidakis and E. Hrushovski, {\it The model theory of difference
	fields}, preprint Universite Paris VII, July 1996.
\medskip

\noindent
[CHP] Z. Chatzidakis, E. Hrushovski, and Y. Peterzil, {\it ${\rm ACFA}_p$},
	in preparation, 1998.

\medskip
\noindent
[Hr]  E. Hrushovski, {\it The Manin-Mumford conjecture and the model
	theory of difference fields}, manuscript, 1996.

\medskip
\noindent
[La] S. Lang, Introduction to Transcendental Numbers, Addison-Wesley, 
	Reading, MA, 1966. 

\medskip
\noindent
[MFvdP] M. Mend\`{e}s France and A. J. van der Poorten, { \it Automata
	and the arithmetic of formal power series}, Acta Arithmetica, 
	{\bf XLVI} (1986) 211 -- 214.

\medskip
\noindent
[Sc]  T. Scanlon, {\it $p$-adic distance from torsion points of 
	semi-abelian varieties}, Journal f\"{u}r die Reine
	und Angewandte Mathematik {\bf 499} (1998), 225 -- 236.

\bigskip
\noindent
MSRI, 1000 Centennial Drive, Berkeley, CA 94720, USA
\smallskip
\noindent
email: {\tt scanlon@msri.org}

\medskip
\noindent
Dept. of Mathematics, Univ. of Texas, Austin, TX 78712, USA
\smallskip
\noindent
e-mail: {\tt voloch@math.utexas.edu}

\end